\documentclass[12pt]{article}

\begin{document}

\makeatletter
\newcommand{\l@abcd}[2]{\hbox{} \hbox to\textwidth{#1\hfil #2}}
\makeatother
\sloppy
\newcommand{\del}{\Delta}
\newcommand{\ot}{\otimes}
\newcommand{\cf}{{\cal F}}
\newcommand{\es}{ e^{\frac{1}{2}\sigma_{1+1}}}
\newcommand{\mes}{ e^{-\frac{1}{2}\sigma_{1+1}}}
\newcommand{\mmes}{ e^{-\sigma_{1+1}}}
\newcommand{\dee}{\Delta_{\epsilon \epsilon J}}
\newcommand{\de}{\Delta_{\epsilon J}}


\begin{titlepage}

\begin{center}
       {\large Ananikian D.N.${\dagger}$\footnote{e-mail:
        ananikia@grad.physics.sunysb.edu},
        Kulish P.P.${\ddagger}$ \footnote{ e-mail: kulish@pdmi.ras.ru},
        Lyakhovsky V.D.${\dagger}$ \footnote{e-mail:
        lyakhovs@pobox.spbu.ru}}\\
        \vskip0.7cm
        ${\dagger}${\it Theoretical Department, St. Petersburg State
        University,}\\
        {\it 198904, St. Petersburg, Russia}\\
        \vskip0.5cm
        ${\ddagger}${\it  St. Petersburg Department of the Steklov
        Mathematical Institute,}\\
        {\it 191011, St.Petersburg, Russia}
\end{center}
\vspace{1cm}

\centerline{\Large \bf Full chains of twists for}
\vspace{0.3cm}
\centerline{\Large \bf symplectic algebras
\footnote{This  work was
supported by Russian Foundation for Basic Research under the
following grants: N 99-01-00101 (K.P.P.)  and N 00-01-00500
(L.V.D). }}
\vspace{1cm}
\begin{abstract}
The problem of constructing the explicit form for full twist
deformations of simple Lie algebras $g$ with twist carriers
containing the maximal nilpotent subalgebra ${\bf N}^+(g)$ is
studied. Our main tool is the sequence of regular subalgebras
$g_i$ in $U(g)$ that become primitive under the action of extended
jordanian twists ${\cal F_E}: U(g)\longrightarrow U_{\cal E}(g)$.
It is demonstrated that the structure of the sequence $\left\{ g_i
\right\}$ is defined by the extended Dynkin diagram of algebra $g$.
To construct the injection $g_i \subset U_{\cal E}(g)$ the special
deformations of algebras $U_{\cal E}(g)$ are performed. They are
reduced to the (cohomologically trivial) twists ${\cal F}_{s}$.
Thus it is proved that full chains of twists can be written in the
canonical form ${\cal F}_{\cal B}
=\prod\limits^{\stackrel{N}{\leftarrow}}_{} {\cal F}_{{\cal
E}'_i} $. The links ${\cal F}_{{\cal E}'_i}$ in such chains must
contain not only the extended twists ${\cal F_E}$ but also the
factors ${\cal F}_{s}$ whose form depend on the type of the series
of classical algebra $g$. The explicit forms of universal ${\cal
R}$-matrices (and the $R$-matrices in the fundamental
representations) corresponding to full chains of twists for
classical simple Lie algebras are found. The properties of the
construction are illustrated by the example of full chain of
extended twists for algebra $sp(3)$.
\end{abstract}
\end{titlepage}

\section{Introduction}
Triangular Hopf algebras ${\cal A} (m,\Delta ,S,\eta
,\epsilon;{\cal R})$ \cite{FRT} play an essential role in quantum
group theory and applications. Quantization of antisymmetric
$r$-matrices, $r=-r_{21}$, (solutions of classical Yang-Baxter
equation (CYBE)) form an important class of such algebras. The
corresponding triangular solutions of quantum Yang-Baxter equation
(QYBE) can be constructed from the classical ones by means of
Campbell-Hausdorff series \cite{DR83}. However, these
constructions are obviously inappropriate for an efficient usage
of quantum ${\cal R}$-matrices. If one provides the initial Lie
algebra $g$ with primitive costructure $\Delta^{\rm prim}$ and
consider the corresponding Hopf algebra $U(g)$, then the solution
${\cal F} \in U(g)^{\otimes 2}$ of the twist equation \cite{DR83}
\begin{equation}
\label{twist-0}
\left( {\cal F}\right) _{12}\left( \Delta^{\rm prim} \otimes {\rm id}
\right) {\cal F}=\left( {\cal F}\right) _{23}\left( {\rm id}\otimes
\Delta^{\rm prim}
\right) {\cal F},
\end{equation}
allows one to find the solution of QYBE, namely:
 ${\cal R_{\cf}}= {\cf_{21}}{\cf}^{-1}$.

For a long time very few types of twists ${\cal F}$ were known in
a closed form (jordanian twist and it's extensions and Reshetikhin
twist) \cite{R, OGI, GIA, KLM}. In particular, the jordanian
twist\cite{OGI}, defined on the Borel subalgebra  $ B(2)=\{H,E\mid
[H,E]=E\} $, has the twisting element
$$
  \Phi _{\cal J}= \exp\{ H\otimes \sigma\},
$$
where $ \sigma = \ln(1+E) $. This twist generates the solution of QYBE:
$ {\cal R}= (\Phi_{\cal J})_{21} \Phi^{-1}_{\cal J} $.
The classical $ r $-matrix $ r=H\land E $ can be associated with it.
It was demonstrated in ~\cite{KLM} that there exist different
extensions of this twist. For the algebra $U(sl(N))$  the
element $ {\cal F}_{\cal E}\in U(sl(N))^{\otimes 2} $ of the form
\begin{equation}
\label{ex-sln}
  {\cal F}_{\cal E}=\Phi_{\cal E} \Phi_{\cal J} =
\left( \prod^{N-2}_{i=2}\Phi_{{\cal E}_{i}}\right) \Phi_{\cal J}
= \exp \left\{ 2\xi\sum_{i=2}^{N-1}{\cal E}_{1,i} \otimes {\cal
E}_{i,N} {e}^{-\widetilde \sigma} \right\} \exp \{
H_{1,N}\otimes\widetilde\sigma \}
\end{equation}
is a solution of (\ref{twist-0}). Here $ H_{1,N}={\cal
E}_{1,1}-{\cal E}_{N,N} $,
$\widetilde\sigma=\frac{1}{2}\ln(1+2\xi E) $ with $ E={\cal E}_{1,N} $
and $ \{{\cal E}_{i,j}\}_{i,j=1,...,N} $ is the standard matrix basis for the
linear algebra. Algebras deformed by extended jordanian twists can
be also obtained via the specific contraction of Drinfeld-Jimbo
quantizations $U_q (sl(N))$ \cite{ABD}.

The minimal subalgebra $g_c \subset g$,
necessary for the twist ${\cal F}_{\cal E}$ to be
defined, is called the {\bf carrier} (sub)algebra for ${\cal F}_{\cal
E}$. The carrier algebra for the extended twist ${\cal F}_{\cal
E}$ is the multidimensional Heisenberg subalgebra of $sl(N)$ with
generators $ \{{\cal E}_{1,i},{\cal E}_{j,N} \mid
i=2,...,N-1,j=1,...,N-1\}$ extended by the Cartan element
$H_{1,N}$.

In the general case the solutions of the twist equation do not
form an algebra with respect to multiplication in ${\cal A}$; the
product of twists must not be a twist. However, under certain
conditions the compositions of extended twists constitute the
solutions of the twist equation, as it was demonstrated in
~\cite{KLO}. For infinite series of simple Lie algebras the twists
(called ``chains") were found. They were constructed as products,
where each factor is itself a twisting element for the initial
algebra (see Section 3). The structure of these chains is
determined by the fundamental symmetry properties of the
corresponding root systems and is common for all simple Lie
algebras. Such chains were called {\bf canonical}.

In the case of chains of extended twists one can assume without
loss of generality that their carrier algebras belong to the Borel
subalgebra ${\bf B}^+(g)= {\bf N}^+ + {\bf H}$. The reason of such
a restriction is that the chain structure is based on the
solvability property of the carrier subalgebra. The Cartan subalgebra
${\bf H}$ is not always entirely contained in $g_c$. To find the
subalgebra ${\bf H} \cap g_c$ in $g$ it is sufficient to consider
the extended Dynkin diagram for $g$ (see section 4). Thus the
dimension criterium for a chain is the ratio of the dimension ${\rm
dim}{\bf N}^+(g)$ and that of the nil-radical for chain's carrier
algebra. The chain is said to be {\bf full} if $g_c
\supset {\bf N}^+(g)$. Among the maximal canonical chains for
classical Lie algebras only the chains for $U(sl(N))$ appear to be
full. In the case of $ U(sp(N)) $
the carrier subalgebra for canonical chain of extended twists
belongs to the regular subalgebra $ sl(N)\subset sp(N) $ (see
~\cite{KLO}). Such chains are called improper for they are totally
defined by the properties of $ sl(N)$-subalgebra. Maximal
canonical chains for orthogonal algebras are proper but not full.

The purpose of this paper is to construct full proper
chains for all classical Lie algebras.

Let $g^{\bot}_{{\lambda}_0}$ be a subalgebra in $g$ with the root
system orthogonal to the {\bf initial root} $\lambda_0$ of jordanian
twist $\Phi_{{\cal J}_0}$ ($\lambda_0$ is the root of the
generator $E_{\lambda_0}$ in the Borel carrier subalgebra of
$\Phi_{{\cal J}_0}$). Let ${\cal F}_{{\cal E}_0}$ be the maximal
(for $g$) extended twist of the form (\ref{ex-sln}) with the
jordanian factor $\Phi_{{\cal J}_0}$. Consider the Cartan
decomposition for the subalgebra $g^{\bot}_{{\lambda}_0}$ in the
form:
$$
\begin{array}{lcl}
g^{\bot}_{{\lambda}_0}&=&{\bf N}^-(g^{\bot}_{{\lambda}_0}) + {\bf
H}(g^{\bot}_{{\lambda}_0}) + {\bf N}^+(g^{\bot}_{{\lambda}_0})=
\\&&
{\bf N}^-(g^{\bot}_{{\lambda}_0}) + {\bf
B}^+(g^{\bot}_{{\lambda}_0}).
\end{array}
$$
The existence of canonical chains of extended twists is based
on the so-called ``matreshka" effect \cite{KLO}. It amounts to the
primitivization of the costructure of the subalgebra $g^{\bot}_{{\lambda}_0}$
twisted by ${\cal F}_{{\cal E}_0}$. In the canonical
chains of twists for $ B_{n} $ and $ D_{n} $ series (constructed
in ~\cite{KLO}) the subalgebras ${\bf N}^+(g^{\bot}_{{\lambda}_0})
\backslash ({\bf N}^+(g^{\bot}_{{\lambda}_0}) \cap g_c)$ are
nontrivial, thus the canonical chain can not be full. The reason
for this peculiarity is that the coproducts in the space of the subalgebra
$g^{\bot}_{{\lambda}_0}$ are nontrivially
deformed by the preceding links of the
extended twists. It was demonstrated in ~\cite{KL} that the
deformed universal enveloping algebras $U(g^{\bot}_{{\lambda}_0})$
contain not only the subalgebras $g^{\bot}_{{\lambda}_0}$ with
deformed costructure but also the primitive ones equivalent to
$g^{\bot}_{{\lambda}_0}$. This property of orthogonal classical
Lie algebras allows to construct full chains of twists for them,
i.e. chains with the carrier algebras $g_c \supset
N^+(g^{\bot}_{{\lambda}_0})$ \cite{Full}.

In Section 4 the decomposition of the root system consistent with
the structure of extended twists for simple Lie algebra will be
constructed. There the system of positive roots $\Lambda^+$ will
be presented as a union of the {\bf initial root} $\lambda_0$, the
{\bf constituent roots} $\left\{ \lambda' , \lambda''| \lambda' +
\lambda''=\lambda_0 \right\}$ for $\lambda_0$ and the
subsystem $\Lambda^{\bot}_{\lambda_0}$ of positive roots
orthogonal to $\lambda_0$
$$
\Lambda^+ =\lambda_0 \bigcup \{ \lambda'  \} \bigcup \{ \lambda''
\} \bigcup \Lambda^{\bot}_{\lambda_0}.
$$
The systematic construction of full chains is based on this
decomposition. For symplectic algebras we investigate the
existence of proper and full chains of twists (see Section 5). 
Our main tool is the deformed carrier space mentioned above. As soon
as full chains of twists can not be based on the sequences of
classical injections $sp(1)\subset sp(2)\subset ...\subset
sp(N-1)\subset sp(N)$ we shall look for the necessary
injections in the universal enveloping algebras $U(sp(N))$. The
problem can be reduced to the construction of a sequences of
injections $ U(sp(1))\subset ...\subset U(sp(N-1))\subset
U(sp(N))$ where to construct the images of a subalgebra  $sp(M)$
means to construct it's nonlinear (in terms of
initial generators) realization in $U(sp(M+1))$.
First we consider the maximal extended jordanian twists $ {\cal
F}_{{\cal E}_{k}} $ for $U(sp(N-k))$. The full proper chain ${\cal
F}_{{\cal B}_{0\prec (N-1)}} $ of such twists for $ U(sp(N)) $
will be realized as the product of factors $ {\cal F}_{{\cal
E}'_{k}} $. Such construction can be considered as a
generalization of the deformed jordanian twists used to define
chains for orthogonal Lie algebras \cite{Full}. The universality
of the primitivization effect for subalgebras
$g^{\bot}_{{\lambda}_0}$ will be established and the recursion
formula for full chains of twists for $U(sp(N))$ will be obtained.

To find the explicit expressions for full chains of twists for all
classical Lie algebras we propose to introduce the additional
twisting factors. Despite their cohomological triviality they
realize the transitions from the deformed to primitive subalgebras
$g^{\bot}_{{\lambda}_0}$ in $U_{ {\cal F}_{{\cal E}'_{k}} }$
(see Section 6). This
allows to include the deformed carrier spaces (specific to the
full chains construction) in the general scheme so that the latter
remains similar to the canonical one.

The expressions for the appropriate multiparametric universal
${\cal R} $-matrices and $R$-matrices in the fundamental
representation are presented. As an example the case of $
U(sp(1))\subset U(sp(2))\subset U(sp(3)) $ is considered in Section 8.

\section{Basic definitions}
In this section we introduce the necessary notations, remind the
definitions and properties of twists. The algebras are considered
below over the field $C$ of complex numbers.

A Hopf algebra ${\cal A} (m,\del,\eta,\epsilon,S)$ \cite{DR86,FRT,DR90}
with multiplication $m: {\cal A}\otimes {\cal A} \rightarrow {\cal A}$,
coproduct $\del : {\cal A} \rightarrow {\cal A} \otimes {\cal A}$,
unit $\eta : C \rightarrow {\cal A}$, counit $\epsilon: {\cal A}
\rightarrow C$ and antipode $S: {\cal A} \rightarrow {\cal A}$ can
be transformed ~\cite{DR90} by an invertible (twisting) element
${\cal F}\in {\cal A}\ot {\cal A}$, $ {\cal F}=\sum f_{i}^{(1)}
\ot f_{i}^{(2)}$ into a twisted one ${\cal A}_{\cal
F}(m,\del_{\cal F},\epsilon, S_{\cal F})$. Hopf algebra ${\cal
A}_{\cal F}$ has the same multiplication, unit and counit but the
twisted coproduct and antipode given by:
\begin{equation}
\label{newcos}
\del_{\cal F}(a)={\cal F}\del (a){\cal F}^{-1}, \quad
S_{\cal F}(a)= vS(a)v^{-1},
\end{equation}
with
$$
 v=\sum f_{i}^{(1)} S (f_{i}^{(2)}), \qquad a\in {\cal A}.
$$
To provide the necessary properties (coassociativity) for the coproduct
$\del_{\cal F}$ it is sufficient that
the twisting element ${\cal F}$ is a solution of the twist
equations \cite{DR83}:
\begin{equation}
\label{twist-eq}
\begin{array}{l}
\left( {\cal F}\right) _{12}\left( \Delta \otimes {\rm id}
\right) {\cal F}=\left( {\cal F}\right) _{23}\left( {\rm id}\otimes \Delta
\right) {\cal F},
\\ \left( \epsilon \otimes {\rm id}\right) {\cal F}=\left(
{\rm id}\otimes \epsilon \right) {\cal F}=1.
\end{array}
\end{equation}
If the initial algebra ${\cal A}$ is quasitriangular with the
universal element ${\cal R}$, then such is the twisted one ${\cal
A}_{\cal F}(m,\del_{\cal F},\eta, \epsilon, S_{\cal F}, {\cal
R}_{\cal F})$ with the ${\cal R}$-matrix
\begin{equation}
\label{newr}
{\cal R_{\cf}}= {\cf_{21}}{\cal R}{\cf}^{-1}.
\end{equation}

Any quantization of a Lie bialgebra with antisymmetric classical
$r$-matrix can be described by a twist  ~\cite{DR83}. Consider a
skew solution of the classical Yang-Baxter equation (CYBE) $r\in
\land^2 {\bf g}$. Let ${\bf l} \subset {\bf g}$ be the minimal Lie
subalgebra such that the bilinear form $b$ on ${\bf l }$ induced
by $r$ is non-degenerate. This form $b$ is a 2-cocycle for the
algebra ${\bf l}$, $b\in {\bf Z}^2({\bf l},{\bf C})$, i.e.
\begin{equation}
\label{form}
b([x,y],z)+
b([z,x],y)+
b([y,z],x)=0.
\end{equation}
Algebra ${\bf l}$ is called {\bf quasi-Frobenius subalgebra} (of the
initial Lie algebra ${\bf g}$), if there exists a non-degenerate
2-cocycle $b$ defined on it. Algebra ${\bf l}$ is called
{\bf Frobenius}, if there exists a linear functional $c\in {\bf l}^*$,
such that the form $c([x,y])=b(x,y)$ is non-degenerate. If ${\bf
l}$ is a Frobenius algebra, then  $b(x,y)\in {\bf B}^2
({\bf l}, C)$, i.e. the form $b$ is  2-coboundary. The classification of
quasi-Frobenius subalgebras for $sl(n)$ is given in
~\cite{stol14}.

In Section 7 we shall demonstrate how the cocycles of
quasi-Frobenius algebras can be directly used in twist
constructions.

\section{Extended twists}
Extended jordanian twists are associated with the parametric set
${\cal L}=\{{\bf L}(\alpha,\beta)_{\alpha + \beta = 1}\}$ of
Frobenius algebras,
\begin{equation}
\label{l-norm}
  \begin{array}{l}
    [H,E] = \delta E,\\
   \left[ A,B\right] = \gamma E,
  \end{array}
  \begin{array}{l}
    [H,A] = \alpha A, \\
   \left[E,A\right] = [E,B] = 0,
 \end{array}
  \begin{array}{l}
    [H,B] = \beta B, \\
    \alpha + \beta = \delta .
  \end{array}
\end{equation}
The equation (\ref{twist-eq}) for the carrier algebra $g_c={\cal
L}$ has the solutions
\begin{equation}
\label{twist1}
{\cal F}_{{\cal E}(\alpha, \beta)}=\exp\{ A\otimes B {\rm
e}^{-\beta\sigma}\}
\exp\{ H\otimes \sigma \}
\end{equation}
and
\begin{equation}
\label{twist2}
{\cal F}_{{\cal E'}(\alpha, \beta)}=\exp\{-B\otimes A {\rm
e}^{-\alpha\sigma}\}
\exp\{ H\otimes \sigma \}.
\end{equation}
The twists ${\cal F}_{{\cal E}(\alpha, \beta)}$ and ${\cal
F}_{{\cal E'}(\alpha, \beta)}$ correspond to the classical
$r$-matrix, $r=H\land E + A\land B$. Algebras of type ${\bf L}$
can be found in any simple Lie algebra $g$ with ${\rm rank}(g)>1$.

It was demonstrated in \cite{KLO}, that (for classical Lie algebras)
the twists of type (\ref{twist1}) or (\ref{twist2}) can be
systematically composed into sequences named {\bf chains of twists}.
This possibility is based on the existence of sequences of regular
injections
\begin{equation}
\label{inject}
g_p\subset g_{p-1}\ldots \subset g_1\subset g_0=g
\end{equation}
and the properties of the invariant symmetric forms on the carrier
spaces for extended twists. To construct a chain it is necessary
to choose the  {\bf initial root}  $\lambda^k_0$ in the root system
$\Lambda(g_k)$. In the root space $V_{\Lambda(g_k)}$ consider
the subspace $V^{\bot}_{{\lambda}^k_{0}}$ orthogonal to the initial root
$\lambda^k_0$. The root system $\Lambda (g_{k+1})=\Lambda
(g_{k})\bigcap V^{\bot}_{{\lambda}^k_0}$ defines the subalgebra
$g_{k+1}$. Consider the set $\pi_k$ of {\bf constituent roots} for
the root $\lambda^k_0$:
\begin{equation}
\label{c-roots}
\begin{array}{l}
\pi_k =\left\{ \lambda ^{\prime },\lambda ^{\prime \prime }\mid \lambda
^{\prime }+\lambda ^{\prime \prime }=\lambda^k_0;\quad \lambda ^{\prime
}+\lambda^k_0,\lambda ^{\prime \prime }+\lambda^k_0\notin \Lambda \left(
g_k\right) \right\} \\
\pi_k =\pi_k ^{\prime }\cup \pi_k ^{\prime \prime };\qquad \pi_k ^{\prime }
=\left\{ \lambda ^{\prime }\right\} ,\pi_k^{\prime \prime }=\left\{ \lambda
^{\prime \prime }\right\} .
\end{array}
\end{equation}
For a classical Lie algebra $g$ in $g^{\bot}_{{\lambda}^k_0}$ one
can always find a subalgebra $g_{k+1}\subseteq
g^{\bot}_{{\lambda}^k_0}\subset g_k$, whose generators have the
primitive costructure after the twist ${\cal F}_{\cal E}$ have
been applied to $U(g_k)$ \cite{KLO}. The phenomenon of
primitivization of a subalgebra $g_{k+1}\subset g_k$ in $U_{\cal
E}(g_k)$ was called the ``matreshka" effect. It allows one to
construct chains of extended jordanian twists,
\begin{equation}
\label{chain-ini}
\begin{array}{l}
{\cal F}_{{\cal B}_{0\prec p}}=\prod_{\lambda ^{\prime }\in \pi _p^{\prime
}}\exp \left\{ E_{\lambda ^{\prime }}\otimes E_{\lambda _0^p-\lambda
^{\prime }}e^{-\frac 12\sigma _{\lambda _0^p}}\right\} \cdot \exp
\{H_{\lambda _0^p}\otimes \sigma _{\lambda _0^p}\}\,\cdot \\ \prod_{\lambda
^{\prime }\in \pi _{p-1}^{\prime }}\exp \left\{ E_{\lambda ^{\prime
}}\otimes E_{\lambda _0^{p-1}-\lambda ^{\prime }}e^{-\frac 12\sigma
_{\lambda _0^{p-1}}}\right\} \cdot \exp \{H_{\lambda _0^{p-1}}\otimes \sigma
_{\lambda _0^{p-1}}\}\,\,\cdot \\
\ldots \\
\prod_{\lambda ^{\prime }\in \pi _0^{\prime }}\exp \left\{ E_{\lambda
^{\prime }}\otimes E_{\lambda _0^0-\lambda ^{\prime }}e^{-\frac 12\sigma
_{\lambda _0^0}}\right\} \cdot \exp \{H_{\lambda _0^0}\otimes \sigma
_{\lambda _0^0}\}\,.
\end{array}
\end{equation}
Here $H_{\lambda _0^k}$ is the generator dual to the root $\lambda _0^k$,
$\sigma _{\lambda_0^k}=\ln(1+E_{\lambda _0^k})$.

Chains of twists give an opportunity to quantize explicitly the
large number of $r$-matrices corresponding to Frobenius
subalgebras in simple Lie algebras \cite{stol14}.

Some peculiarities were found in the construction of chains for
the orthogonal classical algebras. In this case the subalgebra
$g^{\bot}_{{\lambda}^k_0}$ is isomorphic to the direct sum
$sl(2)\oplus g_{k+1}$. After the application of the twist ${\cal
F}_{{\cal E}_k}$ the summand $sl(2)$ obtains the nontrivial
costructure, whereas the coproducts for
$g_{k+1}=so(N-4(k+1))$ remains primitive.
Such primitive summands compose the carrier space for the
canonical chains of type (\ref{chain-ini}) (they are described in
\cite{KLO}). These chains are based on the sequences of injections
$so(N-4(k+1)) \subset so(N-4k)$. It was demonstrated in
\cite{KL, Full} that $U_{{\cal E}_k}(g^{\bot}_{{\lambda}^k_0})$
contains not only the deformed subalgebra  $U_{{\cal
E}_k}(sl^{(k)}(2))$, but also the primitive Hopf subalgebra
isomorphic to $U(sl^{(k)}(2))$. The generators of this algebra can
be composed using the scalar products in the space of
vector subrepresentations $d^v(so(N-4(k+1)))$.
Such subrepresentations appear in the reduction
 of the adjoint representation
of $so(N)$ to the subalgebra $so(N-4(k+1))$.
These subrepresentations have
the weight diagram formed by the projection of constituent roots
(\ref{c-roots}) on the space $V(g^{\bot}_{{\lambda}^k_0})$. As a result the
primitive subalgebra $g^{\bot}_{{\lambda}^k_0}$ can be found in
$U_{{\cal E}_k}(g^{\bot}_{{\lambda}^k_0})$ providing the evidence
that the primitivization of the subalgebra orthogonal to the
initial root is a general property of chain twist deformations. In
this sense the chains of extended twists for $B_N$ series are
similar to that for $A_N$ series. The difference is that in the
$B_N$-case the primitive subalgebra equivalent to
$g^{\bot}_{{\lambda}^k_0}$ is realized on the ``deformed
carrier space" \cite{KL}.

For symplectic simple Lie algebras the situation appears to be
more complicated. When the twist ${\cal F}_{{\cal E}_{k-1}}$ is
applied to $g^{\bot}_{{\lambda}^k_0}\subset U(sp(N-k))$ most of
generators acquire the nonprimitive coproducts. The possibility to
construct the proper chain of twists like (\ref{chain-ini}) fails.
However, as it will be shown below, the primitive subalgebra
isomorphic to $U(sp(N-(k+1)))$ does exist and its generators
can be obtained by
the nonlinear transformations in the space of the subalgebra
$U_{{\cal E}_{k-1}}(sp(N-(k+1)))$.

\section{Extended Dynkin diagrams and chains of twists}

Extended Dynkin diagrams for infinite series of simple Lie
algebras ${\bf g}$ allow one to establish the structure of maximal
chains of twists for the sequence of regular injections and also
define the form of such sequences. Let us remind that the
extended Dynkin diagrams allow one to describe all regular
subalgebras of ${\bf g}$ \cite{VINO}.

{\tolerance=400 Consider the algebra $A_N$ and the corresponding
extended Dynkin diagram: }

{\unitlength=1mm
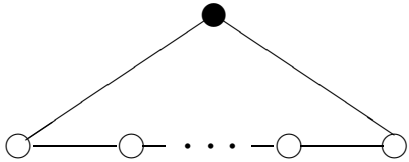
\begin{picture}(50,30)
\put(5,5){\circle{3}}
\put(7,5){\line(1,0){11}}
\put(20,5){\circle{3}}
\put(21.5,5){\line(1,0){3}}
\put(27.5,5){\circle*{0.8}}
\put(30.6,5){\circle*{0.8}}
\put(33.5,5){\circle*{0.8}}
\put(36,5){\line(1,0){3}}
\put(41,5){\circle{3}}
\put(42.5,5){\line(1,0){11}}
\put(55,5){\circle{3}}
\put(31,22.5){\circle*{3}}
\put(31,22.5){\line(-3,-2){25}}
\put(31,22.5){\line(3,-2){24}}
\end{picture}}

Let $ \left\{ \alpha_i  \right\} $ be the set of simple roots for
the algebra $sl(N+1)$, $\theta$ be the grey root: $\theta \equiv
-\sum\limits^N_{i=1}\alpha_i$. In the orthonormal basis $\{
e_i\}_{i=1,...,(N+1)}$ of the space ${\bf R}^{N+1}$ we have
$\alpha_i= e_i-e_{i+1}$ and $-\theta = e_1 -e_{N+1}$. If one
choose the root $-\theta$ to be the initial one, then it follows
from the Dynkin diagram that $g^{\bot}_{\theta}\approx sl(N-1)$,
i.e. the twist ${\cal F}_{{\cal E}_1}$, based on the long root
$-\theta$, will have two constituent basic roots $\alpha_1$ and
$\alpha_N$ in the set $\pi_1$ (see the definition
(\ref{c-roots})). The chain of regular injections for $A_{N}$
series will have the following form:
$$
  sl(2)\subset sl(4)\subset...\subset sl(N-1)\subset sl(N+1)
$$
or
$$
  sl(3)\subset sl(5)\subset...\subset sl(N-1)\subset sl(N+1).
$$
For $N=2k-1$ ($N=2k$) the factor
$\Phi_{{\cal J}_k}
= \exp \{ H_{\lambda _0^k}\otimes \sigma
_{\lambda _0^k}\}$
(correspondingly ${\cal F}_{{\cal E}_k}
=\exp \left\{ E_{\lambda ^{\prime }}\otimes E_{\lambda
_0^k-\lambda ^{\prime }}e^{-\frac 12\sigma _{\lambda
_0^k}}\right\} \cdot \exp \{H_{\lambda _0^k}\otimes \sigma
_{\lambda _0^k}\}$) will be the last in the
maximal chain of twists.

Let us remind that a chain of twists for an algebra $g$ is said to
be full, if it's carrier subalgebra contains the nilpotent
subalgebra $N^{+}(g)$. In the case of $A_N$-series this means that
the maximal canonical chains are full.

Extended Dynkin diagrams for $B_N$ and $D_N$ series have
respectively the following form:

{\unitlength=1mm
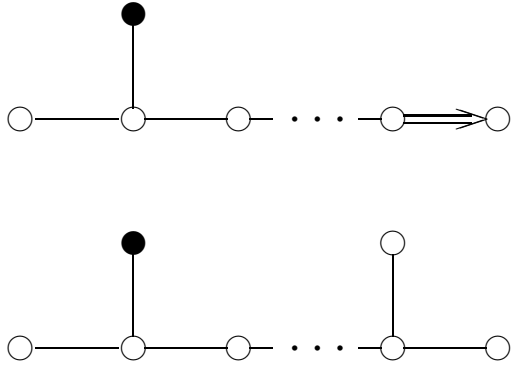
\begin{picture}(220,22)
\put(5,5){\circle{3}}
\put(7,5){\line(1,0){11}}
\put(20,5){\circle{3}}
\put(20,6.5){\line(0,1){11}} \put(21.5,5){\line(1,0){11}}
\put(20,19){\circle*{3}}
\put(34,5){\circle{3}}
\put(35.5,5){\line(1,0){3}}
\put(41.5,5){\circle*{0.8}} \put(44.5,5){\circle*{0.8}}
\put(47.5,5){\circle*{0.8}} \put(50,5){\line(1,0){3}}
\put(54.5,5){\circle{3}}
\put(67,5){\line(-3,1){4}}
\put(67,5){\line(-3,-1){4}} \put(56,5.5){\line(1,0){9}}
\put(56,4.5){\line(1,0){9}} \put(68.5,5){\circle{3}}
\end{picture}
}

{\unitlength=1mm
\begin{picture}(110,25)
\put(5,5){\circle{3}}
\put(7,5){\line(1,0){11}}
\put(20,5){\circle{3}}
\put(20,6.5){\line(0,1){11}} \put(21.5,5){\line(1,0){11}}
\put(20,19){\circle*{3}}
\put(34,5){\circle{3}}
\put(35.5,5){\line(1,0){3}}
\put(41.5,5){\circle*{0.8}} \put(44.5,5){\circle*{0.8}}
\put(47.5,5){\circle*{0.8}} \put(50,5){\line(1,0){3}}
\put(54.5,5){\circle{3}}
\put(56,5){\line(1,0){11}}
\put(54.5,6.5){\line(0,1){11}} \put(68.5,5){\circle{3}}
\put(54.5,19){\circle{3}}
\end{picture}
}

\noindent
It follows immediately that the structure of the spaces orthogonal to the
initial root and the corresponding subalgebras are
$g^{\bot}_{\theta}(B_N)\approx
B_{N-2}\oplus sl(2)$ and $g^{\bot}_{\theta}(D_N)\approx
D_{N-2}\oplus sl(2)$. Thus, the chains of injections for $B_{N}$
and $D_{N}$ series have the form: $$ so(3)\subset
so(7)\subset...\subset so(2N-3)\subset so(2N+1) $$ or $$
so(5)\subset so(9)\subset...\subset so(2N-3)\subset so(2N+1) $$
and
$$
  so(4)\subset so(8)\subset...\subset so(2(N-2))\subset so(2N)
$$
or
$$
  so(6)\subset so(10)\subset...\subset so(2(N-2))\subset so(2N).
$$
For $B_N$ series $\Phi_{{\cal J}_k}$ will be the last factor in
the maximal chain of twists (where $N=2k-1$ if $N$ is odd, and
$N=2(k-1)$ if $N$ is even).

The $D_N$-series is remarkable for the appearance of two
possibilities (depending on $N$) in the last but one step of a
maximal chain of twists:

{\unitlength=1mm
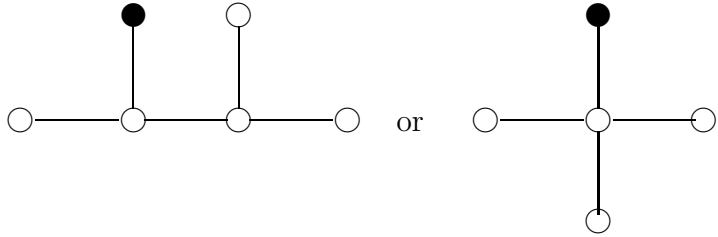
\begin{picture}(60,34)
\put(5,15){\circle{3}}
\put(7,15){\line(1,0){11}}
\put(20,15){\circle{3}}

\put(20,16.5){\line(0,1){11}}
\put(21.5,15){\line(1,0){11}}
\put(20,29){\circle*{3}}
\put(34,15){\circle{3}}
\put(35.5,15){\line(1,0){11}}

\put(48.5,15){\circle{3}}
\put(34,16.5){\line(0,1){11}} \put(34,29){\circle{3}}
\put(55,13.5){or}
\end{picture}
}
{\unitlength=1mm
\begin{picture}(110,34)
\put(5,15){\circle{3}}
\put(7,15){\line(1,0){11}}
\put(20,15){\circle{3}}
\put(20,16.5){\line(0,1){11}}
\put(22,15){\line(1,0){11}}
\put(20,29){\circle*{3}}
\put(20,13.5){\line(0,-1){11}}
\put(20,1.5){\circle{3}}
\put(34,15){\circle{3}}
\end{picture}}

\noindent The first diagram corresponds to even-odd orthogonal
algebra, and the second one to even-even algebra. In the first
case the subalgebra $g^{\bot}_{\theta}$ is isomorphic to
$sl(4)\oplus sl(2)$ and $\Phi_{{\cal J}_k}$ ($N=2k-1$) is the last
factor in the maximal chain. In the second case
$g^{\bot}_{\theta}$ is isomorphic to $sl(2)\oplus sl(2)\oplus
sl(2)$ and $\Phi^1_{{\cal J}_k}$, $\Phi^2_{{\cal J}_k}$,
$\Phi^3_{{\cal J}_k}$ ($N=2k$) are the last three  factors in the
maximal chain of twists. It was shown in \cite{Full}, that for
odd $N$ ($D_{2k+1}$ algebras) in $so(2N)$ there always exists
the independent Cartan generator which does not belong to the
carrier algebra for the maximal chain of twists. In other cases, for so(2N)
with even  $N$ and algebras of series $B$, the carrier subalgebra for the
maximal chain of twists coincides with the Borel subalgebra ${\bf
B}^{+}(g)$.

For $C_N$ series the extended Dynkin diagram has the form

{\unitlength=1mm
\begin{picture}(50,10)
\put(5,5){\circle*{3}}
\put(18.5,5){\line(-3,1){4}}
\put(18.5,5){\line(-3,-1){4}}
\put(6.5,4.5){\line(1,0){9}}
\put(6.5,5.5){\line(1,0){9}}
\put(20,5){\circle{3}}

\put(21.5,5){\line(1,0){3}}
\put(27.5,5){\circle*{0.8}}
\put(30.6,5){\circle*{0.8}}
\put(33.5,5){\circle*{0.8}}
\put(42.5,5){\line(3,1){4}}
\put(42.5,5){\line(3,-1){4}}
\put(36,5){\line(1,0){3}}
\put(41,5){\circle{3}}
\put(44.5,5.5){\line(1,0){9}}
\put(44.5,4.5){\line(1,0){9}}
\put(55,5){\circle{3}}
\end{picture}
}

\noindent Here the subalgebra orthogonal to the initial root is
isomorphic to $C_{N-1}$ and the jordanian twist $\Phi_{{\cal
J}_N}$ is the last factor in the maximal chain of twists. The
Borel subalgebra ${\bf B}^{+}(sp(N))$ appears to be the carrier
subalgebra for the maximal chain.

{\bf Proposition.} The system of positive roots $\Lambda^+$ for
classical Lie algebra $g$ with Cartan decomposition
$$ g = {\bf
N}^-(g) + {\bf H}(g) + {\bf N}^+(g)=  {\bf N}^-(g) + {\bf B}^+(g)
$$
can be presented as a union of the initial root $\lambda_0$,
its constituent roots $\left\{ \lambda' , \lambda''| \lambda' +
\lambda''=\lambda_0 \right\}$ and the subsystem
$\Lambda^{\bot}_{\lambda_0}$ of positive roots orthogonal to
$\lambda_0$:
\begin{equation}
\label{decomp}
\Lambda^+ =\lambda_0 \bigcup \{ \lambda' \}
 \bigcup \{ \lambda'' \}
\bigcup \Lambda^{\bot}_{\lambda_0}.
\end{equation}
$\spadesuit$

{\bf  Proof.} Any highest root  $\lambda_0= -\theta$ is long. If
${\rm rank}(g)>1$ then the basic root $\alpha_i$ can not be collinear to
the highest one. For any positive root $\alpha \in \Lambda^+$ we
have
$$ (\lambda_0 , \alpha) \geq 0.
$$
The case of zero
projection corresponds to the subsystem
$\Lambda^{\bot}_{\lambda_0}$. For any $\left\{\alpha |(\lambda_0 ,
\alpha) > 0 , \alpha \neq
\lambda_0 \right\}$ the difference $\lambda_0 - \alpha $ is a
positive root. Indeed
\begin{equation}
\label{proj} (\lambda_0 - \alpha,\lambda_0)
 = (\lambda_0,\lambda_0)(1- \frac{1}{2}m(\alpha,\lambda_0)),
\end{equation}
where the entry of the Cartan matrix $m(\alpha,\lambda_0)$ in our
situation has a single value equal to 1 \cite{BUR}. Hence, any
positive root $\alpha$ noncollinear and nonorthogonal to
$\lambda_0$ is constituent for $\lambda_0$ and the appropriate
pair is $\left\{\alpha, \lambda_0 - \alpha \right\}$. $\spadesuit$

We have shown that the set $\Lambda^+ \backslash
\Lambda^{\bot}_{\lambda_0}$ consists of the root $\lambda_0$ and
its constituent roots. It is obvious (see (\ref{proj})) that any
constituent root has a projection on $\lambda_0$ equal to $1/2$.

The decomposition (\ref{decomp}) can be considered as an algebraic
foundation for the system of full chains of extended twists for
simple algebras.

\section{Chains of twists for $ U(sp(N)) $}

Consider the sequence of injections
\begin{equation}
\label{inj}
 U(sp(1))\subset...\subset U(sp(N-1))\subset U(sp(N)).
\end{equation}
In the root system
$$
\Lambda( sp(N))  =\{\pm e_{i} \pm e_{j}, \pm 2e_{i}\}; \quad
i,j=1,...,N
$$
choose the vector $\lambda_0 = 2e_{1} $ to be the
initial root. Introduce the following notations:
\begin{enumerate}
\item  $e_{i \pm j} \equiv e_i \pm e_j, \quad i \neq j$\,\,\,  for the short
roots;
\item $ E_{i\pm j} $ and
$ E_{i+i} $ (respectively $ F_{i\pm j} $ and $ F_{i+i} $) --- the
generators corresponding to the roots $ e_{i\pm j} $ and $ 2e_{i}
$ (respectively $ -e_{i\pm j} $ and $ -2e_{i} $), here $ i<j $;
\item $ \pi_{1} $ --- the set of constituent roots for
$\lambda_0$,
$$ \pi_{1}=\{\lambda',\lambda''\mid \lambda' +\lambda''=2e_{1}; \, \,
\lambda'+2e_{1}, \lambda''+2e_{1}\notin \Lambda \}.
$$
\end{enumerate}
Obviously the constituent roots are $
\pi_{1}=\{e_{1-i}, e_{1+i}\mid i=2,\ldots,N\} $.

Let $ H_{ii} $ be the Cartan generators dual to the roots $
2e_{i} $. The Borel subalgebra for $ sp(N) $ in the basis $
\{E_{i+i}, E_{i\pm j}, H_{ii}\mid i,j=1,\ldots,N;\, i<j\} $ is
given by the commutation relations
$$
\begin{array}{l}
\,[H_{ii},E_{n+n}]=\delta _{in}E_{n+n}, \\
\,[H_{ii},E_{n+m}]=\frac 12(\delta _{in}+\delta _{im})E_{n+m}, \\
\,[H_{ii},E_{n-m}]=\frac 12(\delta _{in}-\delta _{im})E_{n-m}, \\
\end{array}
\begin{array}{l}
\,[E_{i-j},E_{n+n}]=2\delta _{jn}E_{i+n}, \\
\,[E_{i-j},E_{n-m}]=\delta _{jn}E_{i-m}-\delta _{mi}E_{n-j}, \\
\,[E_{i-j},E_{n+m}]=\delta _{jn}E_{i+m}+\delta _{jm}E_{i+n}.
\end{array}
$$ It is natural to start the construction of a chain by the
jordanian twist $$
  \Phi _{{\cal J}_{1}}= \exp\{H_{11}\otimes \sigma_{1+1}\}, \quad
  \sigma_{1+1} = \ln(1+E_{1+1}).
$$
Since the triple of roots $ \{ 2e_{1}, e_{1-i},e_{1+i}|
i=2,\ldots, N
\} $ defines the subalgebra $ {\bf L}^{i}(\alpha,\beta) $ with $
\alpha=\beta=\frac{1}{2} $, the constituent roots $ \{
e_{1-i},e_{1+i}| i=2,\ldots, N \} $ allow one to construct the
full extension $\Phi_{{\cal E}_{1}}$ of the twist $\Phi _{{\cal
J}_{1}}$: $$
  \Phi_{{\cal E}_{1}}= \prod^{\stackrel{N}{\leftarrow}}_{i=2}
  e^{E_{1-i}\otimes E_{1+i} e^{-\frac{1}{2}\sigma_{1+1}}}.
$$
The extended twist
$$
 {\cal F}_{{\cal E}_1}=\Phi_{{\cal
E}_{1}}\Phi_{{\cal J}_{1}},
$$
performs the deformation $U(sp(N))
\longrightarrow U_{{\cal E}_1}(sp(N))$. The costructure of the
subalgebra  $ U_{{\cal E}_1}(sp(N-1)) \subset U_{{\cal
E}_1}(sp(N)) $ is defined by the following relations:
\begin{equation}
\label{costr}
\begin{array}{lcl}
\Delta_{{\cal E}_1} (E_{i+i}) & = & E_{i+i} \otimes 1 + 1\otimes E_{i+i}
+2E_{1+i}\otimes E_{1+i} e^{-\frac{1}{2}\sigma_{1+1}}+ \\
&  & + E_{1+1}\otimes E^{2}_{1+i}
e^{-\sigma_{1+1}}, \\
\Delta_{{\cal E}_1} (E_{i+j}) & = & E_{i+j} \otimes 1 + 1\otimes E_{i+j}+
E_{1+i} \otimes E_{1+j} e^{-
\frac{1}{2}\sigma_{1+1}}+ \\  &  & + E_{1+j} \otimes E_{1+i} e^{-
\frac{1}{2}\sigma_{1+1}}+ E_{1+1} \otimes E_{1+i}E_{1+j}
e^{-\sigma_{1+1}},
\\
\Delta_{{\cal E}_1} (E_{i-j}) & = & E_{i-j} \otimes 1 + 1\otimes E_{i-j},
\\
\Delta_{{\cal E}_1} (F_{i+i}) & = & F_{i+i} \otimes 1 + 1\otimes F_{i+i}
+2E_{1-i}\otimes E_{1-i} e^{-
\frac{1}{2}\sigma_{1+1}}- \\  &  & - E_{1-i}^{2}\otimes E_{1+1}
e^{-\sigma_{1+1}}, \\
\Delta_{{\cal E}_1} (F_{i+j}) & = & F_{i+j} \otimes 1 + 1\otimes F_{i+j} +
E_{1-i} \otimes E_{1-j} e^{-
\frac{1}{2}\sigma_{1+1}}+ \\  &  & + E_{1-j} \otimes E_{1-i} e^{-
\frac{1}{2}\sigma_{1+1}} - E_{1-i}E_{1-j}\otimes E_{1+1} e^{-\sigma_{1+1}},
\\ \Delta_{{\cal E}_1} (F_{i-j}) & = & F_{i-j} \otimes 1 + 1\otimes
F_{i-j},
\end{array}
\end{equation}
where $ i,j =2,\ldots,N ; \quad j>i  $. Other generators of the
subalgebra $U_{{\cal E}_1}(sp(N-1))$ retain the former primitive
coproducts.

In order to construct a chain of twists corresponding to the
sequence (\ref{inj}), it is necessary to find the carrier algebra
for the jordanian twist (and it's extensions) containing the
generator of the initial root $\lambda^1_0=2e_2$  for
$sp(N-1)$ subalgebra. It is
generally assumed that the subalgebra to be twisted must possess
the primitive costructure, since in this case one may profit by
using the known solution of the twist equation. The expressions
(\ref{costr}) indicate that in our case it is impossible to apply
the standard procedure of chain construction. The problem will be
solved, if one can find in $ U_{{\cal E}_1}(sp(N)) $ a primitive
subalgebra isomorphic to $ U(sp(N-1)) $ and generated by the
vectors of the space $sp(N-1)$ and the element $\sigma_{1+1}$.

Such a subalgebra does exist in $ U_{{\cal E}_1}(sp(N)) $. To make
sure of this consider the following set of generators:
\begin{equation}
\label{new-gen}
\begin{array}{lcl}
E^{\prime}_{i+i} & = & E_{i+i} - E^{2}_{1+i} e^{-\sigma_{1+1}}, \\
F^{\prime}_{i+i} & = & F_{i+i} - E^{2}_{1-i}, \\
E^{\prime}_{i+j}
& = & E_{i+j} - E_{1+i}E_{1+j} e^{-\sigma_{1+1}}, \\
F^{\prime}_{i+j} & = & F_{i+j} - E_{1-i}E_{1-j}, \\
E^{\prime}_{i-j} & = & E_{i-j},\\ F^{\prime}_{i-j} & = &
F_{i-j},\\ H^{\prime}_{ii} & = & H_{ii}.
\end{array}
\end{equation}
Using the formulas (\ref{costr}) it is not difficult to check that
the costructure of the generators $X^{\prime}$ is primitive. The
space with the basis $ E^{\prime}_{i\pm j}, F^{\prime}_{i\pm j} ,
H^{\prime}_{i+i} \quad (i,j =2,...,N ;\, j
\geq i)$ forms an algebra $sp'(N-1)$ equivalent to $ sp(N-1) $.

The existence  of the primitive algebra $sp'(N-1) \subset
U(sp(N-1))$ verifies the universality of the ``matreshka" effect.
The main point of this effect is that in the universal
enveloping algebra deformed by the full extended twist there
exists the primitive subalgebra, its root system consists of all
the roots belonging to the hyperplane orthogonal to $\lambda_0$.
In order to construct explicitly the chain of twists for regular
injections (\ref{inject}), it is necessary to find the
transformation of the basis leading to the primitive subalgebra like  $
sp'(N-1) $. Note, that for the construction of a chain it is not
important how the generators of  $ sp(N)\backslash sp(N-1)
$ are transformed.

Given a cocommutative Hopf algebra $ U(sp'(N-1)) $, generated by
$ E^{\prime}_{i\pm j}, F^{\prime}_{i\pm j} , H^{\prime}_{i+i}
\quad (i,j =2,...,N ;\, j \geq i)$, we are able to perform the twisting
procedure with $ 2e_{2} $ being the initial root, and so on.

In the $k$-th step we shall have a cocommutative Hopf algebra
$U(sp^{k-1}(N-(k-1)))$ and the possibility to apply the
composition
\begin{equation}
\label{jext-k}
{\cal F}_{{\cal E}_k} = \Phi_{{\cal E}_{k}^{k-1}}
\Phi _{{{\cal J}}_{k}^{k-1}}
\end{equation}
which consists of the jordanian twist $$
  \Phi _{{{\cal J}}_{k}^{k-1}}
  = \exp\{H_{kk}\otimes \sigma^{(k-1)}_{k+k}\}
$$
and it's full extension
$$
\Phi_{{\cal E}_{k}^{k-1}}= \prod_{i>k}
 \exp\left\{{E^{(k-1)}_{k-i} \otimes E^{(k-1)}_{k+i}
e^{-\frac{1}{2}\sigma^{(k-1)}_{k+k}}}\right\}.
$$
Here
$\sigma^{(k-1)}_{k+k} = \ln(1+E^{(k-1)}_{k+k}) $.

Applying the sequence of twists of the type (\ref{jext-k}),
\begin{equation}
\label{can-chain} {\cal F}_{{\cal B}_0 \prec
(k-1)}=\prod^{\stackrel{k}{\leftarrow}}_{j=1} {\cal F}_{{\cal
E}_j},
\end{equation}
we shall obtain the subalgebra $U_{{\cal E}_k}(sp^{k-1}(N-k))$
whose costructure will be defined by the following equalities:
\begin{equation}
\label{costr-k}
\begin{array}{lcl}
\Delta_{{\cal E}_k} (E^{(k-1)}_{i+i}) & = & E^{(k-1)}_{i+i} \otimes 1
+ 1\otimes E^{(k-1)}_{i+i}
+2E^{(k-1)}_{k+i}\otimes E^{(k-1)}_{k+i} e^{-
\frac{1}{2}\sigma^{(k-1)}_{k+k}}+ \\  &  & + E^{(k-1)}_{k+k}\otimes
\left( E^{(k-1)}_{k+i}\right) ^{2}
e^{-\sigma^{(k-1)}_{k+k}}, \\
\Delta_{{\cal E}_k} (E^{(k-1)}_{i+j}) & = & E^{(k-1)}_{i+j} \otimes 1
+ 1\otimes E^{(k-1)}_{i+j}+
E^{(k-1)}_{k+i} \otimes E^{(k-1)}_{k+j} e^{-
\frac{1}{2}\sigma^{(k-1)}_{k+k}}+ \\  &  & + E^{(k-1)}_{k+j}
\otimes E^{(k-1)}_{k+i} e^{-
\frac{1}{2}\sigma^{(k-1)}_{k+k}}+ E^{(k-1)}_{k+k}
\otimes E^{(k-1)}_{k+i}E^{(k-1)}_{k+j}
e^{-\sigma^{(k-1)}_{k+k}},
\\
\Delta_{{\cal E}_k} (E^{(k-1)}_{i-j}) & = & E^{(k-1)}_{i-j} \otimes 1
+ 1\otimes E^{(k-1)}_{i-j},
\\
\Delta_{{\cal E}_k} (F^{(k-1)}_{i+i}) & = & F^{(k-1)}_{i+i} \otimes 1
+ 1\otimes F^{(k-1)}_{i+i}
+2E^{(k-1)}_{k-i}\otimes E^{(k-1)}_{k-i} e^{-
\frac{1}{2}\sigma^{(k-1)}_{k+k}}- \\  &  & -\left( E^{(k-1)}_{k-i}\right)
^{2}\otimes E^{(k-1)}_{k+k}
e^{-\sigma^{(k-1)}_{k+k}}, \\
\Delta_{{\cal E}_k} (F^{(k-1)}_{i+j}) & = & F^{(k-1)}_{i+j}
\otimes 1 + 1\otimes F^{(k-1)}_{i+j} + E^{(k-1)}_{k-i} \otimes
E^{(k-1)}_{k-j} e^{- \frac{1}{2}\sigma^{(k-1)}_{k+k}}+ \\  &  & +
E^{(k-1)}_{k-j} \otimes E^{(k-1)}_{k-i} e^{-
\frac{1}{2}\sigma^{(k-1)}_{k+k}} - E^{(k-1)}_{k-i}E^{(k-1)}_{k-j}
\otimes E^{(k-1)}_{k+k} e^{-\sigma^{(k-1)}_{k+k}},
\\ \Delta_{{\cal E}_k} (F^{(k-1)}_{i-j}) & = & F^{(k-1)}_{i-j}
\otimes 1 + 1\otimes
F^{(k-1)}_{i-j}.
\end{array}
\end{equation}
Thus the result of the sequence of nonlinear transformations is
the reproduction in the $k$-th step of the costructure
(\ref{costr}). This allows to continue the construction of the
chain. In $U_{{\cal E}_k}(sp^{k-1}(N-(k-1)))$ there exists the
primitive subalgebra isomorphic to $U(sp(N-k))$. It can be
obtained by the nonlinear basis transformation
\begin{equation}
\label{k-generators}
\begin{array}{lcl}
  E^{(k)}_{i+i}& = &E^{(k-1)}_{i+i} - \left( E^{(k-1)}_{k+i}\right)^{2}
  e^{-\sigma^{(k-1)}_{k+k}},
\\
  E^{(k)}_{i+j}& = &E^{(k-1)}_{i+j} - E^{(k-1)}_{k+i}E^{(k-1)}_{k+j}
  e^{-\sigma'_{k+k}},
\\
  F^{(k)}_{i+i}& = &F^{(k-1)}_{i+i} - \left( E^{(k-1)}_{k-i}\right)^{2},
\\
  F^{(k)}_{i+j}& = &F^{(k-1)}_{i+j} - E^{(k-1)}_{k-i}E^{(k-1)}_{k-j},
\\
 E^{(k)}_{i-j}& = &E^{(k-1)}_{i-j} ,
\\
F^{(k)}_{i-j}& = &F^{(k-1)}_{i-j} ,
\\ H^{(k)}_{ii}& = &H^{(k-1)}_{ii}, \quad \quad i,j =k+1,...,N;\,
j>i.
\end{array}
\end{equation}

Continuing the specified procedure recursively one can obtain the
explicit expression for the maximal chain of extended twists corresponding 
to the sequence of injections (\ref{inj}),
\begin{equation}
\label{f-genchain}
\begin{array}{rcl}
{\cal F}_{{\cal B}_{0 \prec (N-1)}} & = & e^{H_{NN}\otimes
\sigma^{(N-1)}_{N+N}}\times \\
&  & e^{E^{(N-2)}_{(N-1)-N}\otimes E^{(N-2)}_{(N-1)+N} e^{-
\frac{1}{2}\sigma^{(N-2)}_{(N-1)+(N-1)}}} e^{H_{(N-1)(N-1)}\otimes
\sigma^{(N-2)}_{(N-1)+(N-1)}}\times\\  &  & \vdots \\
&  & e^{E'_{2-N}\otimes E'_{2+N}e^{-\frac
12\sigma^{\prime}_{2+2}}}\dots e^{E'_{2-3}\otimes
E'_{2+3}e^{-\frac 12\sigma^{\prime}_{2+2}}} e^{H_{22}\otimes
\sigma^{\prime}_{2+2}}\times \\ &  & e^{E_{1-N}\otimes
E_{1+N}e^{-\frac 12\sigma _{1+1}}}\dots e^{E_{1-2}\otimes
E_{1+2}e^{-\frac 12\sigma_{1+1}}} e^{H_{11}\otimes \sigma _{1+1}}
.
\end{array}
\end{equation}

\section{Basis transformation and trivial twists}
Now we shall consider the inner automorphisms of the deformed algebras 
$U_{{\cal E}_k}(sp(N-k))$ and demonstrate that the primitive generators of
$U_{{\cal E}_k}(sp(N-k-1))$ constructed
in the previous section by means of basis transformations can be also 
obtained as a result of cohomologically trivial twist.

{\bf Proposition.} In the subalgebra
$sp(N-k-1)  \subset U_{{\cal E}_k}(sp(N-k))$
deformed by the $k$-th step of the chain
${\cal F}_{{\cal B}_0\prec (N-1)}$ the coboundary twist
$$
{\cal F}_{s_k}=\left( s_k^{-1}\otimes s_k^{-1}\right) \Delta _{{\cal E}%
_k}\left( s_k\right)
$$
produces the same primitivization of generators as that obtained by
the basis transformation (\ref{k-generators}).
The element $s_k$ that generates this trivial twist has the
following form:
\begin{equation}
\label{s-k}
\begin{array}{c}
s_k=\exp \left\{\left( \sum\limits_{i=k+1}^N E_{k-i}E_{k+i}\right)
f(E_{k+k})\right\},\qquad f(E_{k+k})=-\frac{\sigma _{k+k}}{2E_{k+k}}.
\end{array}
\end{equation}
$\spadesuit $

The twist ${\cal F}_{s_k}$
must be placed after the factor ${\cal F}_{{\cal E}_k}\in {\cal F}_{{\cal
B}_0\prec (N-1)}$. The simbol $\Delta _{{\cal E}_k}$ in
${\cal F}_{s_k}$ denotes the
(initially primitive) coproduct twisted by the $k$-th link of the chain of
extended twists (in the case under consideration this is the
maximal canonical extended twist for the subalgebra $U(sp(N-k))$).

{\bf Proof.} Let $s$ be the element that induces the inner
automorphism for the algebra $U(sp(N))$ that results in the basis
transformation (\ref
{new-gen}). This means that the generators $Q$ of the deformed algebra 
$U_{{\cal E}_1}(sp(N-1))$ are connected with the generators $Q^{\prime
}$ of the primitive subalgebra $U(sp(N-1))$ by the adjoint
transformation
\begin{equation}
\label{l'-l}Q^{\prime }=sQs^{-1}.
\end{equation}
It is easy to check that (up to a scalar factor) the element $s$
have the form:
$$
s=\exp \left\{\left(
\sum_{i=2}^N E_{1-i}E_{1+i}\right) f(E_{1+1})\right\}.
$$
The function
$f(E_{1+1})$ can be fixed using the relations (\ref{l'-l}). For
example, taking $Q=E_{2+2}$ and $Q^{\prime }=E_{2+2}^{\prime }$
one gets:
\begin{equation}
\label{E'-E}
\begin{array}{c}
sE_{2+2}s^{-1}= {e}^{{\rm ad} \left( \left(
\sum_{i=2}^NE_{1-i}E_{1+i}\right) f(E_{1+1})\right) }E_{2+2} \\
=E_{2+2}-E_{1+2}^2{e}^{-\sigma _{1+1}}.
\end{array}
\end{equation}
Hence $f(E_{1+1})$ must satisfy the equation
$$
{e}^{2E_{1+1}f(E_{1+1})}-1=-E_{1+1}{e}^{-\sigma _{1+1}}=-\left( 1-{e}
^{-\sigma _{1+1}}\right) .
$$
It follows that
$$
f(E_{1+1})=-\frac{\sigma _{1+1}}{2E_{1+1}}.
$$

Consider the generators $X\in U(sp(N-1))$. According to the
properties of the type (\ref{E'-E}) and (\ref{k-generators}) we
have the following
relation:
$$
\begin{array}{c}
\Delta _{
{\cal E}_1}\left( X^{\prime }\right) =\Delta ^{{\rm prim}}\left(
X^{\prime }\right) =X^{\prime }\otimes 1+1\otimes X^{\prime }= \\
=(s\otimes s)\left(
\Delta ^{{\rm prim}}X\right) (s\otimes s)^{-1}.
\end{array}
$$
Thus
$$
\begin{array}{c}
\Delta
{\rm ^{prim}}X=(s\otimes s)^{-1}\Delta _{{\cal E}_1}\left(
X^{\prime }\right) (s\otimes s)= \\ =(s\otimes s)^{-1}\Delta
_{{\cal E}_1}\left( s\right) \Delta _{{\cal E}_1}\left( X\right)
\Delta _{{\cal E}_1}\left( s^{-1} \right) (s\otimes s)\\
={\cal F}_s \Delta _{{\cal E}_1}\left( X\right) ({\cal F}_s )^{-1}.
\end{array}
$$
The expression ${\cal F}_s=(s\otimes s)^{-1}\Delta _{{\cal E}%
_1}\left( s\right) $ is the twisting element of the
cohomologically trivial
twist. If we apply this twist deformation to the subalgebra 
$U_{{\cal E}_1}(sp(N-1))$ 
its initial (nondashed) generators acquire the
primitive coproducts. So when we start to compose the next link of
the chain the necessary twists are functions of
the generators of the
undeformed universal enveloping subalgebra (here it is $U(sp(N-1))$).

The arguments presented above are valid for any link of the chain.
In particular one can avoid the basis transformation
(\ref{k-generators}) in the $k$-th step of the chain ${\cal
F}_{{\cal B}_0\prec (N-1)}$ by inserting the factor ${\cal
F}_{s_k}=\left( s_k^{-1}\otimes s_k^{-1}\right) \Delta _{{\cal
E}_k}\left( s_k\right) $ after (to the left of)  the factor ${\cal
F}_{{\cal E}_k}$. Here $s_k$ is given by the formula (\ref{s-k}).
$\spadesuit $

As a result the maximal chain of twists can be written in the form
\begin{equation}
\label{fullbas}{\cal F}_{{\cal B}_0\prec (N-1)}=
\prod_{i=1}^{\stackrel{N}{\leftarrow }}{\cal F}_{s_i}{\cal F}_{{\cal E}_i}.
\end{equation}
Notice that here the extended twists ${\cal F}_{{\cal E}_k}$ have
the canonical form (\ref{chain-ini}) in terms of the initial generators of $
sp(N)$. This transcription for the maximal chain of twists vividly
demonstrates the universality of the ``matreshka'' effect.

The coboundary twists can be used to compensate
the ``deformations" of the spaces
$V_{\lambda _0}^{\bot }$ in the chains of extended twists for any
classical Lie algebra. For the linear algebras ( $A_N$ series) these
twists
correspond to the identical transformation, ${\cal F}_{s_i}=1\otimes 1$,
and the expression (\ref{fullbas}) describes the canonical
chain.

For orthogonal algebras the coboundary twists are defined by
the elements
$$
s_k=e^{ \left(
E_{2k-1}E_{2k}+\sum_{i=2k+1}^N\left(
E_{(2k-1)+(i)}E_{(2k)-(i)}+E_{(2k-1)-(i)}E_{(2k)+(i)}\right) \right)
f(E_{(2k-1)+(2k)})}
$$
for $ g\in B_N $ and
$$
s_k=e^{  \sum_{i=2k+1}^N\left(
E_{(2k-1)+(i)}E_{(2k)-(i)}+E_{(2k-1)-(i)}E_{(2k)+(i)}\right)
f(E_{(2k-1)+(2k)})}
$$
for $ g\in D_N $.
In both cases ($g\in B_N$ or $D_N$) the twist ${\cal F}_{s_k}$ deforms
only the $sl(2)$-subalgebra in $g_{\lambda _0^k}^{\bot }$. The
latter conserves its primitivity when two twists ${\cal F}_{{\cal E}_k}$ 
and ${\cal F}_{s_k}$ are successively applied.

Thus the expression (\ref{fullbas}) (with the appropriate choice
of $s_k$) describes the maximal chain of twists for an arbitrary
simple Lie algebra of the infinite series.


\section{$R$-matrices and forms}
The coboundary twists ${\cal F}_{{s}_k}$
do not contribute to classical $r$-matrices. The
expression (\ref{f-genchain}) allows one to use the formula
(\ref{newr}) while finding the explicit form of the quantum
${\cal R}$-matrix for the deformed algebra $U_{{\cal B}_{0\prec
(N-1)}}(sp(N))$:
\begin{equation}
\label{quantr}
  {\cal R}_{{\cal B}_{0\prec (N-1)}}=
\left( {\cal F}_{{\cal B}_{0\prec (N-1)}}\right)_{21}
\left( {\cal F}_{{\cal B}_{0\prec (N-1)}}\right)^{-1}.
\end{equation}
To fix the classical limit it is necessary to introduce
deformation parameters. This can be performed by
scaling the generators of the twisted
algebra $U_{{\cal B}_{0\prec (N-1)}}(sp(N))$:
$$
E_{i+k}\longrightarrow \xi\eta_{k}
E_{i+k};\quad k\ge i.
$$
Considering the expression (\ref{quantr})
in the small neighborhood of $\xi=0$ one can extract the classical
$r$-matrix
\begin{equation}
\label{klassr}
r=\sum_{k=1}^{N}\eta_{k}(H_{kk}\land E_{k+k}+ \kappa_{k}\sum_{i=k+1}^{N}
E_{k-i}\land E_{k+i}).
\end{equation}
Here the multipliers $\eta_{k}$ are arbitrary. They
signify the independence of the chain links in the quasiclassical
limit. The discrete parameter $\kappa _k$ (with values $\kappa
_k=0,1$) indicates whether the extension of the $k$-th jordanian
twist has been included in the chain. (The degenerate ``matreshka"
effect takes place also for the jordanian twist without
extensions.) The number of discrete parameters $\kappa _k$ is
equal to  $N-1$. For full chains of twists we have $\kappa
_k=1$ for all $k$'s.

The Frobenius forms (they are 2-coboundaries)
\begin{equation}
\label{gen-form}
\omega _p=\sum\limits^p_{k=1}\gamma_k E^{*}_{k+k}([\cdot,\cdot]),
\qquad p=1,...,N,
\end{equation}
correspond to the solution (\ref{can-chain}) of the twist equation.
These forms are (by the definition) nondegenerate on the carrier algebra
of the twist. The matrix $r=\omega^{-1}_p$ is defined on the same subalgebra.

For $\eta_k = \frac{1}{\gamma_k}$ in (\ref{klassr}) the
solution of CYBE corresponding to the matrix $r=\omega^{-1}_N$ coincides
with (\ref{klassr}).

Algebra ${\bf B}^+(sp(N))$ has a nontrivial second cohomology
group $H^2({\bf B}^+(sp(N)))$. It is not difficult to be
verify that $H^*_{ii}\land H^*_{kk}$ is a nontrivial 2-cocycle. Therefore
the forms
\begin{equation}
\label{omegahat}
\hat \omega _p=\omega_p + \xi_{ik} H^*_{ii}\land H^*_{kk},
\quad \quad    i,k\leq N,
\end{equation}
are the nondegenerate 2-cocycles. We can associate twists to
the forms $\hat \omega _p$. To find them notice that
application of the chain ${\cal F}_{{\cal B}_{0\prec
(p-1)}}$ gives rise to the set of primitive commuting elements $\{
\sigma^{(k-1)}_{k+k}\mid k=1,...,p \}$. It follows that the algebra
$U_{{\cal B}_{0\prec (p-1)}}(sp(N))$ admits the twist
$$
  \Phi _{\cal R}=\exp \{\xi _{ij} \sigma_i \otimes \sigma_j\}, \quad
\sigma _i \in \{ \sigma^{(k-1)}_{k+k}\mid k=1,...,p \}.
$$
Thus the composition $\Phi_{\cal R} {\cal F}_{{\cal B}_{0\prec
(p-1)}}$ satisfies the twist equation for $U(sp(N))$. The
corresponding ${\cal R}$-matrix looks like $$
  {\cal R}_{{\cal R B}_{0\prec (p-1)}}=
\left( \Phi_{\cal R}{\cal F}_{{\cal B}_{0\prec (p-1)}}\right)_{21}
\left( {\cal F}_{{\cal B}_{0\prec (p-1)}}\right)^{-1} \left(
\Phi_{\cal R} \right)^{-1}.
$$
Tending the deformation parameter to zero one can obtain the
classical $r$-matrix
\begin{equation}
\label{hatr}
  \hat r= r + \sum^p_{i=1}\sum^p_ {k=i+1} \zeta_{ik} E_{i+i} \land E_{k+k}.
\end{equation}
If we put $p=N$ in (\ref{omegahat}), the matrix
$$
{\hat \omega_N}^{-1}=r+\sum^N_{i=1}\sum^N_{k=i+1}
\frac{\xi_{ik}}{\gamma_i \gamma_k} E_{i+i}\land E_{k+k}
$$
coincides (up to redifinition of parameters) with the $r$-matrix (\ref{hatr}).

As it was demonstrated in Section 6 devoted to cohomologically
trivial twists for any simple
Lie algebra of four infinite series
there exists the general expression for the twisting
element of the full chain. To construct the appropriate quantum
system one needs the $R$-matrix in a finite-dimensional representation.
In particular, one can use the fundamental
representations of lowest dimension --- the
first fundamental representations. Notice that in these
representations the basis transformations that we need are degenerate,
i.e. the dashed generators coincide with the initial
ones. This is the direct consequence of the generators nilpotency in the
first fundamental representation.

Let
$$
d(E_{\lambda})=L_{\lambda}, \qquad
d(H_{\lambda})=M_{\lambda}
$$
be the generators of the Borel subalgebra in the first fundamental
representation $d$. Then for any series of Lie
algebras the $R$-matrix in this representation will look like:
\begin{equation}
\label{fund}
\begin{array}{lcl}
d^{\otimes 2}(\cal R)& = &
\prod\limits_{\lambda'\in\pi'_{k}}^{\rightarrow}
e^{L_{\lambda_0^k-\lambda'}\otimes L_{\lambda'}}
e^{L_{\lambda_0^k}\otimes M_{\lambda_0^k}}\dots
\prod\limits_{\lambda'\in\pi'_{0}}^{\rightarrow}
e^{L_{\lambda_0-\lambda'}\otimes L_{\lambda'}}
e^{L_{\lambda_0}\otimes M_{\lambda_0}}\\
& & \prod\limits_{\lambda'\in\pi'_0}^{\leftarrow}
e^{ - M_{\lambda_0} \otimes L_{\lambda_0} }
e^{ - L_{\lambda'} \otimes L_{\lambda_0^k-\lambda'} }\dots
\prod\limits_{\lambda'\in\pi'_k}^{\leftarrow}
e^{ -M_{\lambda_0^k} \otimes L_{ \lambda_0^k } }
e^{-L_{\lambda'}\otimes L_{\lambda_0^k-\lambda'}}.
\end{array}
\end{equation}
Here it was taken into account that in this
representation we have
$d(\sigma_{\lambda_0^k})=\ln(1+L_{\lambda_0^k})=L_{\lambda_0^k} $.

Let us introduce the deformation parameters by the substitution:
$$
  L_{\lambda_0^k}\rightarrow \xi\eta_k L_{\lambda_0^k},
  \qquad
  L_{\lambda_0^k-\lambda'}\rightarrow \xi\eta_k L_{\lambda_0^k-\lambda'}.
$$
(The other generators remain unchanged and it is easy to
check that this is the algebra
automorphism.) Due to the fact that the only nonzero terms in the
jordanian twists and the extensions in the representation $d^{\otimes 2}$
are linear in $\xi$ we get the simple expression for $d^{\otimes 2}(\cal R)$:
\begin{equation}
\begin{array}{lcl}
\label{fun1}
d^{\otimes 2}(\cal R) & = & [
(1+\xi\eta_k\sum\limits_{\lambda'\in\pi'_k}L_{\lambda_0^k-\lambda'}
\otimes L_{\lambda'})(1+\xi\eta_k L_{\lambda_0^k}\otimes M_{\lambda_0^k})
]\times\\
& &\vdots\\
& &
[
(1+\xi\eta_0\sum\limits_{\lambda'\in\pi'_0}L_{\lambda_0-\lambda'}
\otimes L_{\lambda'})(1+\xi\eta_0 L_{\lambda_0}\otimes M_{\lambda_0})
]\times\\
& &
[
(1-\xi\eta_k M_{\lambda_0}\otimes L_{\lambda_0})(1-\xi\eta_0\sum\limits_
{\lambda'\in\pi'_0}L_{\lambda'}\otimes L_{\lambda_0-\lambda'})
]\times\\
& &\vdots\\
& &
[
(1-\xi\eta_0 M_{\lambda_0^k}\otimes L_{\lambda_0^k})(1-\xi\eta_k\sum\limits_
{\lambda'\in\pi'_k}L_{\lambda'}\otimes L_{\lambda_0^k-\lambda'})
].
\end{array}
\end{equation}
\newcommand{\lmk}{L_{\lambda^k_0-\lambda'_k}}
\newcommand{\lmi}{L_{\lambda^i_0-\lambda'_i}}
\newcommand{\lmj}{L_{\lambda^j_0-\lambda'_j}}
\newcommand{\lk}{L_{\lambda'_k}}
\newcommand{\lj}{L_{\lambda'_j}}
\newcommand{\li}{L_{\lambda'_i}}
\newcommand{\la}{\lambda}
{\bf Lemma}:
\begin{equation}
\begin{array}{l}
\label{fund2}
d^{\otimes 2}({\cal R})
 =  I \otimes I+
\xi\sum\limits_{k=0}^{N}\eta_k\{\sum\limits_{\lambda'\in\pi'_k}
(L_{\lambda_0^k-\lambda'}\land L_{\lambda'_k})
+L_{\lambda_0^k}\land M_{\lambda_0^k}\}+\\
+ \xi^2 \sum\limits_{k}\eta_k^2\{\sum\limits_{\la'_k\in\pi_k}
(\lmk\lk\otimes\lk\lmk)+L_{\la_0^k}M_{\la_0^k}
\otimes M_{\la_0^k} L_{\la_0^k} \}+\\
 +  \xi^2 \sum\limits_{i<j,\la'_i\in\pi_i}\eta_i\eta_j
(L_{\la^j_0}\li\otimes M_{\la_0^j}\lmi + \lmi M_{\la_0^j}
\otimes\li L_{\la_0^j})+\\
 +  \xi^2\sum\limits_{i<j,\la'_i\in\pi_i,\la'_j\in\pi_j }\eta_i\eta_j
(\lmj\li \otimes \lj\lmi + \lmi\lj\otimes\li\lmj).
\end{array}
\end{equation}
$\spadesuit$

The proof can be carried out by induction through the length of the chain.

Apply the Lemma to the full chains of extended twists for the
algebras $A_{N-1}$. Replacing all the generators by the appropriate
$N\times N$ basic matrices, one can see that only the
terms originating from the product of jordanian twists give the nonzero
contribution of the order $\xi^2$:
\begin{equation}
\begin{array}{l}
d^{\otimes 2}({\cal R})= I\ot I +\\

+ \xi \sum\limits_{k=0}^{[\frac{N}{2}]}\eta_k \{
\sum\limits_{i=2}^{N-k}({\cal E}_{i,N-k}\land {\cal E}_{k+1,i}) +
{\cal E}_{k+1,N-k}
\land ({\cal E}_{k+1,k+1} - {\cal E}_{N-k,N-k} ) \} +\\
+ \xi^2 \sum\limits_{k=0}^{[\frac{N}{2}]}\eta_k^2 {\cal
E}_{k+1,N-k} \ot {\cal E}_{k+1,N-k}.
\end{array}
\end{equation}
For symplectic series all the three groups of terms give the nontrivial
contribution of the second order:
\begin{equation}
\begin{array}{l}
d^{\otimes 2} ({\cal R})= I\otimes I + \xi \{
\sum\limits_{i=1}^N \eta_i {\cal E}_{i,(N+i)}\land ({\cal E}_{i,i}-
{\cal E}_{(N+i),(N+i)})+\\
+\sum\limits_{j=2}^{N}\sum\limits_{i<j}\eta_{j}
({\cal E}_{i,(j+N)}+{\cal E}_{j,(i+N)})\land ({\cal E}_{i,j}-
{\cal E}_{(j+N),(i+N)} ) \} + \\
+\xi^2 \sum\limits_{k=1}^{N}(N-k+1)\eta_k^2 {\cal
E}_{k,\frac{N}{2}+k}
\ot {\cal E}_{k,\frac{N}{2}+k}+\\
+\xi^2 \sum\limits_{i<j}^{N}\eta_i \eta_j
({\cal E}_{i,N+j}\ot {\cal E}_{i,N+j} + {\cal E}_{j,N+i}\ot {\cal
E}_{j,N+i})+\\
+\xi^2 \sum\limits_{i<j}^N (N-j) \eta_i \eta_j
({\cal E}_{i,N+j}\ot {\cal E}_{i,N+j} + {\cal E}_{j,N+i}\ot {\cal
E}_{j,N+i}).
\end{array}
\end{equation}
Here, for convenience the indices of $\eta$ are shifted by
one, i.e. the substitution $\eta_k\rightarrow \eta_{k+1}$ is performed.


\section{Chains of twists for $sp(3)$}
As an example consider the full chain of twists for the algebra
$sp(3)$. This is the simplest case where the specific
structure of symplectic series and it's difference from
linear and orthogonal algebras can be seen. Consider the root system
$$
\Lambda \left( sp(3)\right)=\{\pm e_i\pm e_j,\quad \pm 2e_i\}
\quad (i,j=1,2,3).
$$
For the initial root  $\lambda _0=2e_1$ the constituent roots are
$\lambda^{\prime }=e_1-e_i$ and $\lambda ^{\prime \prime
}=e_1+e_i$.

The full proper chain of twists for $U(sp(3))$ should be based on
the sequence of regular injections
$$
U(sp(1))\subset U(sp(2))\subset U(sp(3)).
$$
The roots of $\Lambda(sp(2))$ are orthogonal to the long root
$\lambda_{0}\in \Lambda(sp(3))\backslash \Lambda (sp(2))$. In
other words for the appropriate indexation
of basic vectors in $\Lambda$ we have $\Lambda(sp(2))\,\bot\, 2e_{1}$ and
$\Lambda(sp(1))\,\bot \, 2e_{1}, 2e_{2}$.

We start the construction of the full chain of twists by
performing the jordanian twist with the carrier subalgebra
generated by $\left\{H_{11}, E_{1+1} \right\}$:
$$
\Phi _{{\cal J}_1}=\exp\{H_{11}\otimes \sigma _{1+1}\},\qquad
\sigma =\ln(1+E_{1+1}).
$$
The sets $\{2e_1,e_{1-2},e_{1+2}\}$ and
$\{2e_1,e_{1-3},e_{1+3}\}$ define two extensions ${\cal E}^{I}$
and ${\cal E}^{II}$ for $\Phi _{{\cal J}_1}$. So the full
jordanian extended twist has the twisting element
\begin{equation}
\label{exttw} {\cal F}_{{\cal E}_1}=\underbrace{e^{E_{1-3} \otimes
E_{1+3}e^{-\frac 12\sigma _{1+1}}}}_{{\cal E}^{II}}\underbrace{%
e^{E_{1-2}\otimes E_{1+2}e^{-\frac 12\sigma _{1+1}}}}_{{\cal E}%
^{I}}\cdot \underbrace{e^{H_{11}\otimes \sigma _{1+1}}}_{\cal J}.
\end{equation}
Being applied to $U(sp(3))$ this twist deforms the subalgebra
$g_{\lambda^0_0}^{\perp }=g_{(2e_1)}^{\perp }= U_{{\cal E}_{1}}(sp(2))$
which obtains the costructure defined by the relations:
\begin{equation}
\label{costrsp}
\begin{array}{lcl}
\Delta_{{\cal E}_1} (E_{2+2}) & = & E_{2+2} \otimes 1 + 1\otimes E_{2+2}
+2E_{1+2}\otimes E_{1+2} e^{-
\frac{1}{2}\sigma_{1+1}}+ \\  &  & + E_{1+1}\otimes E^{2}_{1+2}
e^{-\sigma_{1+1}}, \\
\Delta_{{\cal E}_1} (E_{2+3}) & = & E_{2+3} \otimes 1 + 1\otimes E_{2+3}+
E_{1+2} \otimes E_{1+3} e^{-
\frac{1}{2}\sigma_{1+1}}+ \\  &  & + E_{1+3} \otimes E_{1+2} e^{-
\frac{1}{2}\sigma_{1+1}}+ E_{1+1} \otimes E_{1+2}E_{1+3} e^{-\sigma_{1+1}},
\\ \Delta_{{\cal E}_1} (E_{2-3}) & = & E_{2-3} \otimes 1 + 1\otimes E_{2-3},
\\
\Delta_{{\cal E}_1} (E_{3+3}) & = & E_{3+3} \otimes 1 + 1\otimes E_{3+3}
+2E_{1+3}\otimes E_{1+3} e^{-
\frac{1}{2}\sigma_{1+1}}+ \\  &  & + E_{1+1}\otimes E^{2}_{1+3}
e^{-\sigma_{1+1}}, \\
\Delta_{{\cal E}_1} (F_{2+2}) & = & F_{2+2} \otimes 1 + 1\otimes F_{2+2}
+2E_{1-2}\otimes E_{1-2} e^{-
\frac{1}{2}\sigma_{1+1}}- \\  &  & - E_{1-2}^{2}\otimes E_{1+1}
e^{-\sigma_{1+1}}, \\
\Delta_{{\cal E}_1} (F_{2+3}) & = & F_{2+3} \otimes 1 + 1\otimes F_{2+3} +
E_{1-2} \otimes E_{1-3} e^{-
\frac{1}{2}\sigma_{1+1}}+ \\  &  & + E_{1-3} \otimes E_{1-2} e^{-
\frac{1}{2}\sigma_{1+1}} - E_{1-2}E_{1-3}\otimes E_{1+1} e^{-\sigma_{1+1}},
\\ \Delta_{{\cal E}_1} (F_{2-3}) & = & F_{2-3} \otimes 1 + 1\otimes F_{2-3},
\\
\Delta_{{\cal E}_1} (F_{3+3}) & = & F_{3+3} \otimes 1 + 1\otimes F_{3+3}
+2E_{1-3}\otimes E_{1-3} e^{-
\frac{1}{2}\sigma_{1+1}}- \\  &  & - E_{1+3}^{2}\otimes E_{1+1}
e^{-\sigma_{1+1}}.
\end{array}
\end{equation}
Now we shall use formulas (\ref{new-gen}) to construct the deformed
generators in the twisted algebra $U_{{\cal E}_1}(sp(2))$,
\begin{equation}
\label{sp-prime}
\begin{array}{lcl}
E^{\prime}_{i+i} & = & E_{i+i} - E^{2}_{1+i} e^{-\sigma_{1+1}}, \\
F^{\prime}_{i+i} & = & F_{i+i} - E^{2}_{1-i}, \\
E^{\prime}_{2+3} & = & E_{2+3} - E_{1+2}E_{1+3} e^{-\sigma_{1+1}}, \\
F^{\prime}_{2+3} & = & F_{2+3} - E_{1-2}E_{1-3}, \\
H^{\prime}_{ii} & = & H_{ii}.
\end{array}
\qquad i=2,3
\end{equation}
The dashed generators are primitive and form the subalgebra
$sp^{\prime}(2)$.

To continue the construction of the chain consider the algebra
$U(sp^{\prime}(2))$ on the deformed carrier space generated by
$E'$, $F'$ and $H'$. In the standard form we introduce for
$sp^{\prime}(2)$ an independent root system $\Lambda(sp^{\prime}(2)) $.

For the algebra $U(sp^{\prime}(2)) $ there exists the jordanian
twist based on the long root $\lambda_{2+2}=2e_{2}\in
\Lambda(sp^{\prime}(2)) $,
\begin{equation}
\label{2-jord}
\Phi_{{\cal J}'_{2}}=\exp\{H_{22}\otimes\sigma^{\prime}_{2+2}\},
\qquad
\sigma^{\prime}_{2+2}= \ln (1+E^{\prime}_{2+2}),
\end{equation}
and it's extension $\Phi_{{\cal E}'_{2}}$ defined by the roots
$\left\{ e_{2-3},\ e_{2+3} \right\}$. In other words the algebra
$U(sp^{\prime}(2)) $ admits the extended jordanian twist with the
element
$$
{\cal F}_{{\cal E}'_2}=\underbrace
{e^{E^{\prime}_{2-3}\otimes E^{\prime}_{2+3} e^{-
\frac{1}{2}\sigma^{\prime}_{2+2}}}}_{{\cal E}'_2} e^{H_{22}\otimes
\sigma^{\prime}_{2+2}}.
$$
Being applied to $U(sp^{\prime}(2)) $
it nontrivially transforms the costructure of the subalgebra
$sp^{\prime}(1) \subset sp^{\prime}(2)$. Let us again pass to the
deformed space $sp^{\prime \prime}(1) \subset U_{{\cal E}_2{\cal
E}_1}(sp(3))$ with primitive generators (compare with the corresponding
expressions in \cite{KL}):
$$
\begin{array}{lcl}
E^{\prime\prime}_{3+3} & = & E^{\prime}_{3+3} - (E^{\prime}_{2+3})^{2}

e^{-\sigma^{\prime}_{2+2}}, \\
F^{\prime\prime}_{3+3} & = & F^{\prime}_{3+3} - (E^{\prime}_{2-3})^{2}, \\
H^{\prime\prime}_{33} & = & H_{33}. \\
&  &
\end{array}
$$

The next link, the last one for the full chain in case of
$U(sp(3))$, is degenerate. It does not contain extensions. The
subalgebra $U(sp''((1))$ on the twice deformed space is generated
by the elements $E^{\prime\prime}$, $F^{\prime\prime}$ and
$H^{\prime\prime} $. It admits the jordanian twist $$
\Phi_{{\cal J}''_{3}}=\exp\{ H_{33}\otimes \sigma^{\prime\prime}_{3+3}\},
\qquad
\sigma^{\prime\prime}_{3+3}= \ln(1+E^{\prime \prime}_{3+3}).
$$
Thus, the jordanian twist
$$
{\cal F}_{{\cal E}''_3}=\Phi_{{\cal J}''_{3}}.
$$
is the last factor in our chain.

The full chain of twists for $U(sp(3))$ will have the following form:
\begin{equation}
\label{f-chain}
\begin{array}{lcl}
{\cal F}_{{\cal B}_{0 \prec 2 }} & = & {e}^{H_{33}\otimes
\sigma^{\prime\prime}_{3+3}}\times \\
&  & \times {e}^{E^{\prime}_{2-3}\otimes E^{\prime}_{2+3} e^{-
\frac{1}{2}\sigma^{\prime}_{2+2}}} {e}^{H_{22}\otimes
\sigma^{\prime}_{2+2}}\times\\
&  & \times{e}^{E_{1-3}\otimes
E_{1+3}{e}^{-\frac 12\sigma _{1+1}}} {e}^{E_{1-2}\otimes E_{1+2}{e}^{-\frac
12\sigma_{1+1}}} {e}^{H_{11}\otimes \sigma _{1+1}}. \\
\end{array}
\end{equation}
Here the generators of  $sp'(2)$ are defined by the formulas
(\ref{sp-prime}) and the generator $ E^{\prime\prime}_{3+3} $
has the form
$$
E^{\prime\prime}_{3+3}=E_{3+3}-E_{1+3}^{2}
e^{-\sigma_{1+1}}-(E_{2+3}-E_{1+2}E_{1+3} e^{-\sigma_{1+1}})^{2}
e^{-\sigma^{\prime}_{2+2}}.
$$

Notice, that the initial sequence of regular injections has also been 
deformed:
$$
U(sp''(1)) \subset U(sp'(2)) \subset U(sp(3)).
$$
Nevertheless this new sequence still consists of  symplectic
regular subalgebras.  Thus the chain (\ref{f-chain}) is
 proper for the algebra $U(sp(3))$ (the improper chains of
twists were described in \cite{KLO}). Algebra ${\bf B}^+(sp(3))$ is
the carrier for the twist (\ref{f-chain}). Consequently this
chain of extended twists is  full.

The quantum ${\cal R}$-matrix, corresponding to the twist
(\ref{f-chain}), and it's classical counterpart are as follows:
\begin{equation}
\label{rmat}
{\cal R}=({\cal F}_{{\cal B}_{0 \prec 2}})_{21} ({\cal F}_{{\cal
B}_{0 \prec 2}})^{-1},
\end{equation}
$$
r=\sum_{k=1}^{3}\eta_{k}(H_{kk}\land E_{k+k}+ \kappa_{k}\sum_{i=k+1}^{3}
E_{k-i}\land E_{k+i}).
$$
The approach based on the cohomologically trivial twists, developed in
the section 6, allows one to rewrite the expression (\ref{f-chain}) in
the form
\begin{equation}
\label{fullbassp}
{\cal F}_{{\cal B}_{0 \prec 2}}
=\prod\limits_{i=1}^{\stackrel{3}{\leftarrow}} {\cal F}_{s_i}{\cal
F}_{{\cal E}_i} ,
\end{equation}
here
$$ s_1 = \exp \{ \left( E_{1-2} E_{1+2}+E_{1-3}
E_{1+3}\right) f(E_{1+1})\},
$$
$$
s_2=\exp \{ E_{2-3} E_{2+3} f(E_{2+2})\},
$$
$$
s_3=1.
$$

Finally let us notice that the recursive procedure proposed
in Section 5 allows one to construct full chains of twists for
any classical Lie algebra.

\section{Conclusions }
Explicitly constructed maximal proper chain of twists for
$U(sp(N))$ terminates the process of composition of full proper
chains for classical Lie algebras. Thus it was demonstrated that
the ``matreshka" effect is universal. The following proposition
can be formulated: the Hopf algebra $U_{\cal E}(g)$ with simple $g$,
deformed by the maximal (canonical) extended twist ${\cal F}_{\cal E}$,
contains a primitive subalgebra equivalent to the universal
enveloping algebra $U(g^{\bot}_{{\lambda}_0})$ where the root
system of $g^{\bot}_{{\lambda}_0}
\subset g$ is orthogonal to the initial root ${\lambda}_0$.
For classical algebras this statement was proved in Section 6.
It would be interesting to check it's validity for the five exceptional
Lie algebras.

There exists a variety of applications where the explicit
expressions for twisting elements received above are essential. The
explicit form of the universal ${\cal R}$-matrix allows one
to get the $R$-matrix in an arbitrary representation. Thus the
first fundamental representation $d$ was exploited in Section 7 and the
vector representation $d_v$ in \cite{Full}.

The twist deformation of a coalgebra induces the
modification of the Clebsh-Gordan coefficients in the decomposition of
the tensor product of irreducible representation $d_{i}\otimes d_{j}$. New
Clebsh-Gordan coefficients are defined by the action of the
matrix operator $F=d_{i}\otimes d_{j} (\cal F)$ of the twisting
element ${\cal F}$ on the initial ones \cite{KulStol}.

Due to the embedding of simple Lie algebras in the appropriate
Yangians $U(g)\subset Y(g)$, there exists a possibility to deform
the Yangians by the chains of twists for the algebras $U(g)$
\cite{KST}.

If we consider the Yangian corresponding to the symplectic algebra
$sp(N)$ the $R$-matrix in the defining representation $d \otimes d
\subset {\rm Mat}(2N,C)\otimes {\rm Mat}(2N,C)$ looks as follows
\cite{RESH2}
$$
 uI + {\cal P} - u(u+N+1)^{-1}\widetilde{\cal K}.
$$
Here $u$ is the spectral parameter, $\widetilde{\cal K}$ is the rank one
projector, that in the basis $({\bf C}^2 \otimes{\bf C}^N)
\otimes ({\bf C}^2 \otimes{\bf C}^N)$ can be written as $P_-{\cal K}$
with $P_-$ being the antisymmetrizer in  ${\bf C}^2  \otimes {\bf
C}^2$ and the operator ${\cal K}$ can be obtained from the permutation 
${\cal P}$ in ${\bf C}^N  \otimes {\bf C}^N$ by transposing the first
tensor multiplier. The twisted $R$-matrix will have the form $$
R_{\cal F} (u)= uF_{21}F^{-1} + {\cal P} - u(u+N+1)^{-1}F_{ 21}
\widetilde{\cal K}F^{-1}. $$ This leads to the deformation of the
Hamiltonian density in the integrable $sp(N)$-spin system (see
\cite{KulStol}),
$$
h_F =F {\cal P}F^{-1} - (N+1)^{-1}F
\widetilde{\cal K} F^{-1} =
      FhF^{-1}.
$$


\end{document}